\documentclass[12pt,reqno]{amsart}

\newtheorem{thm}{Theorem}
\newtheorem{lemma}[thm]{Lemma}

\newcommand{\N}{{\mathbb N}}
\newcommand{\ind}{{\mathbf 1}}

\newcounter{mycount}
\newenvironment{mylist}{\begin{list}{(\roman{mycount})}%
{\usecounter{mycount}\itemsep 0pt}}{\end{list}}

\title{Partition Identities and the Coin Exchange Problem}
\author{Alexander E. Holroyd}
\address{Department of Mathematics, University of British Columbia,
121-1984 Mathematics Rd., Vancouver, BC V6T 1Z2, Canada}
\thanks{Funded in part by an NSERC Discovery Grant and by Microsoft
Research} \email{holroyd(at)math.ubc.ca} 

\begin{document}

\begin{abstract}
The number of partitions of $n$ into parts divisible by $a$ or $b$
equals the number of partitions of $n$ in which each part and each
difference of two parts is expressible as a non-negative integer
combination of $a$ or $b$. This generalizes identities of MacMahon
and Andrews.  The analogous identities for three or more integers
(in place of $a,b$) hold in certain cases.
\end{abstract}

\maketitle

\section{Introduction}
A {\bf partition} of $n$ is an unordered multiset of positive
integers (called {\bf parts}) whose sum is $n$.  For positive
integers $a_1,\ldots,a_m$ we denote the set of non-negative integer
combinations
$$S=S(a_1,\ldots,a_m):=\big\{\!\textstyle\sum_{i=1}^m x_i a_i:
x_1,\ldots ,x_m\in\N\big\},$$
where $\N=\{0,1,\ldots\}$.

\begin{thm}\label{main} For positive integers $n$, $a_1$ and $a_2$,
the following are all equinumerous:
\begin{mylist}
\item partitions of $n$ in which each part and each difference
between two parts lies in $S(a_1,a_2)$; \item partitions of $n$ in
which each part appears with multiplicity lying in $S(a_1,a_2)$;
\item partitions of $n$ in which each part is divisible by $a_1$
or $a_2$.
\end{mylist}
\end{thm}

For example, when $(n,a_1,a_2)=(13,3,4)$, the three sets of
partitions are: (i) $\{(13),(10,3),(7,3,3)\}$; (ii)
$\{(3,3,3,1,1,1,1)\},(2,2,2,1,\ldots,1),$ $(1,\ldots,1)\}$; (iii)
$\{(9,4),(6,4,3),(4,3,3,3)\}$.

We also establish the following partial extension to three or more
integers $a_1,\ldots,a_m$. Let $\sqcap$ and $\sqcup$ denote greatest
common divisor and least common multiple respectively.

\newpage
\begin{thm}
\label{ext} For any positive integers $n$ and $a_1,\ldots,a_m$,
the following are equinumerous:
\begin{mylist}
\item partitions of $n$ in which each part and each difference
between two parts lies in $S(a_1,\ldots,a_m)$; \item partitions of
$n$ in which each part appears with multiplicity lying in
$S(a_1,\ldots,a_m)$.
\end{mylist}
If $a_1,\ldots,a_m$ satisfy
\begin{equation}\tag{*}
\label{star} \forall i=2,\ldots,m,\; \exists j<i \text{ such that
} (a_1\sqcap\cdots\sqcap a_{i-1})\sqcup a_i = a_j\sqcup a_i.
\end{equation}
then in addition the following are equinumerous with (i) and (ii):
\begin{mylist}
\item[(iii)] partitions of $n$ in which each part is divisible by
some $a_i$.
\end{mylist}
\end{thm}

Note that \eqref{star} holds automatically when $m=2$, so Theorem
\ref{main} is a special case of Theorem \ref{ext}.

\section{Remarks}
To avoid uninteresting cases, $a_1,\ldots,a_m$ should be coprime,
and none should be a multiple of another.  (Indeed, if the
greatest common divisor is $g>1$ then Theorem \ref{ext} reduces
easily to the case $(n',a_1',\ldots,a_m')=g^{-1}(n,a_1,\ldots
,a_m)$, while if $a_j$ is a multiple of $a_i$ then the statements
of the theorem are unchanged by removing $a_j$ from
$a_1,\ldots,a_m$).

The set $S$ is sometimes interpreted as describing sums of money
that can be formed using coins of given denominations.  When
$a_1,\ldots,a_m$ are coprime, the complement $S^C:=\N\setminus S$ is
finite; see e.g.\ \cite{ramirez}.  The case $m=2$ was studied by
Sylvester \cite{sylvester}, who proved for $a_1,a_2$ coprime that
$|S^C|=\tfrac12 (a_1-1)(a_2-1)$ and $\max S^C=a_1 a_2-a_1-a_2$. The
case $m\geq 3$ was proposed by Frobenius, and is much less well
understood in general.  An exception is when is when
$a_1,\ldots,a_m$ satisfy a certain condition which is implied by our
condition \eqref{star}; see \cite{nij-wilf}. For more information
see \cite{ramirez}.

\sloppy When $m=2$ we have for example $S(2,3)^C=\{1\}$;
$S(3,4)^C=\{1,2,5\}$; $S(2,5)^C=\{1,3\}$; $S(3,5)^C=\{1,2,4,7\}$;
 $S(4,5)^C$ \linebreak $=\{1,2,3,6,7,11\}$.  Larger sets
$\{a_1,\ldots,a_m\}$ satisfying condition \eqref{star} include
$\{4,6,9\}$; $\{6,8,9\}$; $\{6,9,10\}$; $\{8,12,18,27\}$;
$\{30,42,70,105\}$. We have for instance
$S(4,6,9)^C=\{1,2,3,5,7,11\}$.\fussy

In the case $\{a_1,a_2\}=\{2,3\}$, the equality between (i) and
(iii) in Theorem \ref{main} gives following partition identity due
to MacMahon \cite[\S 299--300]{macmahon}.
\begin{quote}\em
The number of partitions of $n$ into parts not congruent to $\pm 1$
modulo $6$ equals the number of partitions of $n$ with no
consecutive integers and no ones as parts.
\end{quote}
The generalization to $\{a_1,a_2\}=\{2,2r+1\}$, $r\in\N$ was proved
(in a form similar to that above) by Andrews \cite{andrews-gen}.
Other recent work related to MacMahon's identity appears in
\cite{aepr,pnas,hlr}. Somewhat similar identities are proved in
\cite{andrews-lewis}. For more information on partitions and
partition identities see e.g.\ \cite{andrews}.

Finally we note that the second assertion in Theorem \ref{ext}
cannot hold for arbitrary $a_1,\ldots,a_m$ with $m\geq 3$.  For
example, it does not hold for $\{a_1,a_2,a_3\}=\{2,3,5\}$: we have
$S(2,3,5)=S(2,3)$, but allowing multiples of 5 in addition to
multiples of 2 and 3 can clearly increase the number of partitions
of type (iii).

\section{Proofs}

As remarked above, Theorem \ref{main} is the $m=2$ case of Theorem
\ref{ext}.  We will prove the two assertions of Theorem \ref{ext}
separately.  The proofs are simpler when $m=2$, and the reader may
find it helpful to bear this case in mind throughout.

\begin{proof}[Proof of Theorem \ref{ext} (first equality)]
Fix $a_1,\ldots,a_m$, and let $F_n$ and $M_n$ be the sets of
partitions in (i) and (ii) respectively.  We will show that
$|F_n|=|M_n|$.

For a partition $\lambda=(\lambda_1,\ldots,\lambda_r)$ (where
$n=\sum_i\lambda_i$ and $\lambda_1\geq\cdots\geq\lambda_r$), the
{\bf conjugate} partition
$\lambda'=(\lambda'_1,\ldots,\lambda'_{r'})$ is defined as usual
by $r'=\lambda_1$ and $\lambda'_i=\max\{j:\lambda_j\geq i\}$.
Since the set $S$ is closed under addition, the condition that
$\lambda$ has all parts and differences between parts in $S$ is
equivalent to the condition that each {\em adjacent} pair in the
sequence $\lambda_1,\lambda_2,\ldots,\lambda_r,0$ differs by an
element of $S$.  On the other hand, it is readily seen that the
latter condition is equivalent to the condition that $\lambda'$
has all multiplicities in $S$ (indeed this holds for any set $S$).
Hence conjugation is a bijection between $F_n$ and $M_n$.
\end{proof}

Our proof of the second assertion in Theorem \ref{ext} relies on the
two simple lemmas below.  Given integers $a_1,\ldots,a_m$ we write
$$\ell_i:=(a_1\sqcap\cdots\sqcap a_{i-1})\sqcup a_i.$$

\begin{lemma}
\label{gen} If $a_1\ldots,a_m$ satisfy condition \eqref{star} then
we have the formal power series identity
$$\sum_{k\in S(a_1,\ldots,a_m)}q^k = \frac{\prod_{i=2}^m(1-q^{\ell_i})}
{\prod_{i=i}^m(1-q^{a_i})}.$$
\end{lemma}
In the case when $m=2$ and $a_1,a_2$ are coprime, the above
expression has the appealing form $(1-q^{a_1 a_2})
(1-q^{a_1})^{-1}(1-q^{a_2})^{-1}$, as noted in \cite{s-wormald}.
Expressions for the left side for $m=3$ and arbitrary $a_1,a_2,a_3$,
are derived in \cite{denham,s-wormald}.
\begin{proof}[Proof of Lemma \ref{gen}]
We use induction on $m$.  When $m=1$ we have
$$\sum_{k\in
S(a_1)}q^k=1+q^{a_1}+q^{2a_1}+\cdots=\frac{1}{1-q^{a_1}}$$ as
required.

For $m\geq 2$, clearly any $k\in S(a_1,\ldots,a_m)$ can be
expressed as
\begin{equation}\label{weak}
k=xa_m+y, \qquad\text{where }x\in\N \text{ and }y\in
S(a_1,\ldots,a_{m-1}).
\end{equation}
We claim that under condition \eqref{star}, each such $k$ has a
{\em unique} such representation subject to the additional
constraint
\begin{equation}\label{strong}
x<\ell_m/a_m.
\end{equation}
Once this is proved we obtain
$$\sum_{k\in S(a_1,\ldots,a_m)}q^k =
(1+q^{a_m}+q^{2a_m}+\cdots+q^{\ell_m-a_m}) \sum_{k\in
S(a_1,\ldots,a_{m-1})}q^k.$$ By the inductive hypothesis this
equals
$$\frac{1-q^{\ell_m}}{1-q^{a_m}} \times
\frac{\prod_{i=2}^{m-1}(1-q^{\ell_i})}
{\prod_{i=i}^{m-1}(1-q^{a_i})},$$ which is the required
expression.

To check the above claim, let $j=j(m)$ be as in condition
\eqref{star}, and write $d=a_1\sqcap\cdots\sqcap a_{m-1}$, so that
$\ell_m=d\sqcup a_m=a_j\sqcup a_m$. Now note that any representation
$k=xa_m+y$ as in \eqref{weak} that violates \eqref{strong} may be
re-expressed as $k=(x-\ell_m/a_m)a_m+(y+\ell_m)$, where
$x-\ell_m/a_m\in\N$, and $y+\ell_m\in S(a_1,\ldots,a_{m-1})$ (since
$\ell_m$ is a multiple of $a_j$).  By repeatedly applying this we
can reduce $x$ until \eqref{strong} is satisfied, as required.  To
check uniqueness, note that all elements of $S(a_1,\ldots,a_{m-1})$
are divisible by $d$, while the $\ell_m/a_m$ quantities
$0,a_m,2a_m,\ldots,\ell_m-a_m$ are all distinct modulo $d$ (since
$\ell_m=d\sqcup a_m$).  Hence we see that no two distinct
expressions $xa_m+y$ satisfying \eqref{weak},\eqref{strong} can be
equal.
\end{proof}

Let $\ind[\cdot]$ denote an indicator function and let $|$ denote
``divides''.
\begin{lemma}
\label{ind} If $a_1\ldots,a_m$ satisfy condition \eqref{star} then
for any positive integer $k$,
$$\ind\big[a_i|k \text{\rm\ for some }i\big]=\sum_{i=1}^m \ind[a_i|k]
-\sum_{i=2}^m \ind[\ell_i|k].$$
\end{lemma}
When $m=2$ and $a_1,a_2$ are coprime, the lemma is the familiar
inclusion/exclusion formula $\ind[a_1|k \text{ or
}a_2|k]=\ind[a_1|k]+\ind[a_2|k]-\ind[a_1a_2|k]$.
\begin{proof}[Proof of Lemma \ref{ind}]
We use induction on $m$.  The case $m=1$ is trivial.  For $m\geq
2$ we have
\begin{align*}
\ind\big[a_i|k \text{ for some }i\big]=& \ind[a_m|k]
+\ind\big[a_i|k \text{ for some }i<m\big]\\
&-\ind\big[a_m|k, \text{ and }a_i|k \text{ for some }i<m\big]
\end{align*}
We claim that the last condition ``$a_m|k, \text{ and }a_i|k
\text{ for some }i<m$" is equivalent to $\ell_m|k$.  Once this is
established, the result follows by substituting the inductive
hypothesis and the claim into the above equation.

Turning to the proof of the claim, if the given condition holds then
$a_m|k$ and $d|k$, where $d=a_1\sqcap\cdots\sqcap a_{m-1}$.  So $k$
is divisible by $a_m\sqcup d=\ell_m$.  For the converse, recall from
\eqref{star} that $\ell_m=a_m\sqcup a_j$ for some $j<m$, so
$\ell_m|k$ implies $a_m|k$ and $a_j|k$.
\end{proof}

\begin{proof}[Proof of Theorem \ref{ext} (second equality)]
Suppose \eqref{star} holds, and let $M_n$ and $D_n$ denote the
sets of partitions in (ii) and (iii) respectively.  We will show
$|M_n|=|D_n|$.

Using Lemma \ref{gen}, the generating function for $|M_n|$ is
$$G(q):=\sum_{n=0}^\infty |M_n|\, q^n=\prod_{t=1}^\infty
\bigg[ \sum_{k\in S} q^{kt} \bigg]= \prod_{t=1}^\infty
\frac{\prod_{i=2}^{m}(1-q^{\ell_i t})} {\prod_{i=i}^{m}(1-q^{a_i
t})}.
$$
When the product over $t$ is expanded, the factor $(1-q^{\ell_i t})$
contributes a factor $(1-q^k)$ in the numerator for each $k$ that is
a non-negative multiple of $\ell_i$; similarly for the factors in
the denominator. Thus
\begin{align*}
G(q)&=\prod_{k=1}^\infty \Big(1-q^k\Big)^{\textstyle-\sum_{i=1}^m
\ind[a_i|k] +\sum_{i=2}^m \ind[\ell_i|k]}\\
&=\prod_{k=1}^\infty \Big(1-q^k\Big)^{\textstyle-\ind[a_i|k \text{
for some }i]} \;=\!\prod_{ \substack{k\in{\mathbb Z}^+:
\\a_i|k \text{ for some }i}}\frac{1}{1-q^k}.
\end{align*}
(In the second equality we have used Lemma \ref{ind}.) But the
last expression is the generating function for $|D_n|$.
\end{proof}

\section*{Questions}

Can Theorems \ref{main} and \ref{ext} be given simple bijective
proofs?  Dan Romik has found an affirmative answer for Theorem
\ref{main} (personal communication).  Is condition \eqref{star}
necessary and sufficient for the identity between (i) and (iii) in
Theorem \ref{ext}?  For those $a_1,\ldots,a_m$ not satisfying this
identity, are the partitions of type (i) or type (iii) equinumerous
with partitions in some other natural classes? Can condition
\eqref{star} be expressed in a more natural form?

\section*{Acknowledgments}

I thank Dan Romik and George Andrews for valuable comments.

\bibliographystyle{abbrv}
\bibliography{part}

\end{document}